\input amstex
\input amsppt.sty
\magnification=\magstep1
\hsize=30truecc
\vsize=22.2truecm
\baselineskip=16truept
\NoBlackBoxes
\TagsOnRight \pageno=1 \nologo
\def\Z{\Bbb Z}
\def\N{\Bbb N}

\def\C{\Bbb C}
\def\l{\left}
\def\r{\right}
\def\bg{\bigg}
\def\({\bg(}
\def\[{\bg\lfloor}
\def\){\bg)}
\def\]{\bg\rfloor}
\def\t{\text}
\def\f{\frac}

\def\bi{\binom}
\def\eq{\equiv}

\def\ls{\leqslant}
\def\gs{\geqslant}
\def\mo{\roman{mod}}

\def\ve{\varepsilon}

\def\Proof{\noindent{\it Proof}}

\def\Remark{\medskip\noindent{\it  Remark}}

\def\Ack{\medskip\noindent {\bf Acknowledgment}}
\hbox {J. Number Theory 132(2012), no.\,11, 2673--2699.}
\bigskip
\topmatter
\title On sums of Ap\'ery polynomials and related congruences \endtitle
\author Zhi-Wei Sun\endauthor
\leftheadtext{Zhi-Wei Sun}
\rightheadtext{Sums of Ap\'ery polynomials and related congruences}
\affil Department of Mathematics, Nanjing University\\
 Nanjing 210093, People's Republic of China
  \\  zwsun\@nju.edu.cn
  \\ {\tt http://math.nju.edu.cn/$\sim$zwsun}
\endaffil
\abstract The Ap\'ery polynomials are given by
$$A_n(x)=\sum_{k=0}^n\bi nk^2\bi{n+k}k^2x^k\ \ (n=0,1,2,\ldots).$$
(Those $A_n=A_n(1)$ are Ap\'ery numbers.)
Let $p$ be an odd prime. We show that
$$\sum_{k=0}^{p-1}(-1)^kA_k(x)\eq\sum_{k=0}^{p-1}\f{\bi{2k}k^3}{16^k}x^k\pmod{p^2},$$
and that
$$\sum_{k=0}^{p-1}A_k(x)\eq\l(\f xp\r)\sum_{k=0}^{p-1}\f{\bi{4k}{k,k,k,k}}{(256x)^k}\pmod{p}$$
for any $p$-adic integer $x\not\eq0\pmod p$.
This enables us to determine explicitly $\sum_{k=0}^{p-1}(\pm1)^kA_k$ mod $p$,
and $\sum_{k=0}^{p-1}(-1)^kA_k$ mod $p^2$ in the case $p\eq2\pmod3$.
Another consequence states that
$$\sum_{k=0}^{p-1}(-1)^kA_k(-2)\eq\cases 4x^2-2p\pmod{p^2}&\t{if}\  p=x^2+4y^2\ (x,y\in\Z),
\\0\pmod{p^2}&\t{if}\ p\eq3\pmod4.
\endcases$$
We also prove that for any prime $p>3$ we have
$$\sum_{k=0}^{p-1}(2k+1)A_k\eq p+\f 76p^4B_{p-3}\pmod{p^5}$$
where $B_0,B_1,B_2,\ldots$ are Bernoulli numbers.
\endabstract
\thanks 2010 {\it Mathematics Subject Classification}.\,Primary 11A07, 11B65;
Secondary  05A10, 11B68, 11E25.
\newline\indent Supported by the National Natural Science
Foundation (grant 11171140) of China and the PAPD of Jiangsu Higher Education Institutions.
\endthanks
\endtopmatter
\document

\heading{1. Introduction}\endheading

The well-known Ap\'ery numbers given by
$$A_n=\sum_{k=0}^n\bi
nk^2\bi{n+k}k^2=\sum_{k=0}^n\bi{n+k}{2k}^2\bi{2k}k^2\
(n\in\N=\{0,1,2,\ldots\}),$$
play a central role in Ap\'ery's
proof of the irrationality of $\zeta(3)=\sum_{n=1}^\infty 1/n^3$ (see Ap\'ery [Ap] and van der Poorten [Po]).
They also have close connections to modular forms (cf. Ono [O, pp.198--203]).
The Dedekind eta function in the theory of modular forms is defined by
$$\eta(\tau)=q^{1/24}\prod_{n=1}^\infty(1-q^n)\ \quad\t{with}\ q=e^{2\pi i\tau},$$
 where $\tau\in\Bbb H=\{z\in\C:\ \roman{Im}(z)>0\}$ and hence $|q|<1$.
 In 1987 Beukers [B] conjectured that
$$A_{(p-1)/2}\eq a(p)\ (\mo\ p^2)\quad \ \t{for any prime}\ p>3,$$
where $a(n)\ (n=1,2,3,\ldots)$ are given by
$$\eta^4(2\tau)\eta^4(4\tau)=q\prod_{n=1}^\infty(1-q^{2n})^4(1-q^{4n})^4=\sum_{n=1}^\infty a(n)q^n.$$
This was finally confirmed by Ahlgren and Ono [AO] in 2000.

We define Ap\'ery polynomials by
$$A_n(x)=\sum_{k=0}^n\bi nk^2\bi{n+k}k^2x^k=\sum_{k=0}^n\bi{n+k}{2k}^2\bi{2k}k^2x^k\ \ (n\in\N).\tag1.1$$
Clearly $A_n(1)=A_n$. Motivated by the Ap\'ery polynomials,
we also introduce a new kind of polynomials:
$$W_n(x):=\sum_{k=0}^n\bi nk^2\bi{n-k}k^2x^k=\sum_{k=0}^{\lfloor n/2\rfloor}\bi n{2k}^2\bi{2k}k^2x^k\ \ (n\in\N).\tag1.2$$

Recall that Bernoulli numbers $B_0,B_1,B_2,\ldots$ are rational numbers given by
$$B_0=1\ \ \t{and}\ \ \sum^n_{k=0}\bi {n+1}k B_k=0\ \ \ \t{for}\ n\in\Z^+=\{1,2,3,\ldots\}.$$
It is well known that $B_{2n+1}=0$ for all $n\in\Z^+$ and
$$\f x{e^x-1}=\sum_{n=0}^\infty B_n\f{x^{n}}{n!}\ \ \l(|x|<2\pi\r).$$
Also, Euler numbers $E_0,E_1,E_2,\ldots$ are integers defined by
$$E_0=1\ \ \t{and}\ \ \sum^n\Sb k=0\\2\mid k\endSb \bi nk E_{n-k}=0\ \ \ \t{for}\ n\in\Z^+.$$
It is well known that $E_{2n+1}=0$ for all $n\in\N$ and
$$\sec x=\sum_{n=0}^\infty(-1)^n E_{2n}\f{x^{2n}}{(2n)!}\ \ \l(|x|<\f{\pi}2\r).$$

Now we state our first theorem.
\proclaim{Theorem 1.1}
{\rm (i)} Let $p$ be an odd prime. Then
$$\sum_{k=0}^{p-1}(-1)^kA_k(x)\eq\sum_{k=0}^{p-1}(-1)^kW_k(-x)\eq\sum_{k=0}^{p-1}\f{\bi{2k}k^3}{16^k}x^k\pmod{p^2}.\tag1.3$$
Also, for any $p$-adic integer $x\not\eq0\pmod p$, we have
$$\aligned\sum_{k=0}^{p-1}A_k(x)\eq&\sum_{k=0}^{p-1}W_k(x)\pmod{p^2}
\\\eq&\l(\f xp\r)\sum_{k=0}^{p-1}\f{\bi{4k}{k,k,k,k}}{(256x)^k}\pmod{p},
\endaligned\tag1.4$$
where $(-)$ denotes the Legendre symbol.

{\rm (ii)} For any positive integer $n$ we have
$$\f1n\sum_{k=0}^{n-1}(2k+1)A_k(x)=\sum_{k=0}^{n-1}\bi{n-1}k\bi{n+k}k\bi{n+k}{2k+1}\bi{2k}kx^k.\tag1.5$$
If $p>3$ is a prime, then
$$\sum_{k=0}^{p-1}(2k+1)A_k\eq p+\f 76p^4B_{p-3}\  (\mo\ p^5)\tag1.6$$
and
$$\sum_{k=0}^{p-1}(2k+1)A_k(-1)\eq \l(\f{-1}p\r)p-p^3E_{p-3}\pmod{p^4}.\tag1.7$$

{\rm (iii)} Given $\ve\in\{\pm1\}$ and $m\in\Z^+$, for any prime $p$ we have
$$\sum_{k=0}^{p-1}(2k+1)\ve^k A_k^m\eq0\ (\mo\ p).$$
\endproclaim
\Remark\ 1.1. (i) Let $p$ be an odd prime. The author [Su1, Su2] had conjectures on $\sum_{k=0}^{p-1}\bi{2k}k^3/m^k$ mod $p^2$
with $m=1,-8,16,-64,256,-512,4096$.
Motivated by the author's conjectures on $\sum_{k=0}^{p-1}A_k(x)$ mod $p^2$
with $x=1,-4,9$ in an initial version of this paper, Guo and Zeng [GZ, Theorem 1.3] recently showed that
$$\sum_{k=0}^{p-1}A_k(x)\eq\sum_{k=0}^{(p-1)/2}\bi{p+2k}{4k+1}\bi{2k}k^2x^k\pmod{p^2}.$$

(ii) The values of
$$s_n=\f1n\sum_{k=0}^{n-1}(2k+1)A_k\in\Z$$ with $n=1,\ldots,8$ are
$$1,\ 8,\ 127,\ 2624,\ 61501,\ 1552760,\ 41186755,\ 1131614720$$
respectively.
On June 6, 2011 Richard Penner informed the author an interesting application of (1.5):
(1.5) with $x=1$ implies that
$s_n$ is the trace of the inverse of $nH_n$ where $H_n$ refers to the Hilbert matrix $(\f1{i+j-1})_{1\ls i,j\ls n}$.
\medskip

Can we find integers $a_0,a_1,a_2,\ldots$ such that $\sum_{k=0}^{p-1}a_k\eq 4x^2-2p\pmod{p^2}$
if $p=x^2+y^2$ is a prime with $x$ odd and $y$ even? The following corollary provides an affirmative answer!

\proclaim{Corollary 1.1} Let $p$ be any odd prime. Then
$$\aligned&\sum_{k=0}^{p-1}(-1)^kA_k(-2)\eq\sum_{k=0}^{p-1}(-1)^kA_k\l(\f 14\r)
\\\eq&\cases 4x^2-2p\ (\mo\ p^2)&\t{if}\ p\eq1\ (\mo\ 4)\ \&\ p=x^2+y^2\ (2\nmid x),
\\0\ (\mo\ p^2)&\t{if}\ p\eq3\ (\mo\ 4).\endcases
\endaligned\tag1.8$$
\endproclaim
\Proof. It is known (cf. Ishikawa [I]) that
$$\sum_{k=0}^{p-1}\f{\bi{2k}k^3}{64^k}\eq\cases 4x^2-2p\ (\mo\ p^2)&\t{if}\ p\eq1\ (\mo\ 4)\ \&\ p=x^2+y^2\ (2\nmid x),
\\0\ (\mo\ p^2)&\t{if}\ p\eq3\ (\mo\ 4).\endcases$$
The author conjectured that we can replace $64^k$ by $(-8)^k$ in the congruence,
 and this was recently confirmed by Z. H. Sun [S3]. So, applying (1.3) with $x=-2,1/4$
 we obtain (1.8). \qed

\proclaim{Corollary 1.2} Let $p$ be an odd prime. Then
$$\sum_{k=0}^{p-1}A_k\eq c(p)\pmod{p}\tag1.9$$
where
$$c(p):=\cases 4x^2-2p&\t{if}\ p\eq1,3\pmod8\ \&\ p=x^2+2y^2\ (x,y\in\Z),
\\0&\t{if}\ (\f{-2}p)=-1,\ \t{i.e.,}\ p\eq5,7\pmod8.\endcases$$
Also,
$$\aligned&\sum_{k=0}^{p-1}(-1)^kA_k\eq\l(\f{-1}p\r)\sum_{k=0}^{p-1}(-1)^kA_k\l(\f 1{16}\r)
\\\eq&\cases4x^2-2p\ (\mo\ p)&\t{if}\ p\eq1\ (\mo\ 3)\ \t{and}\ p=x^2+3y^2\ (x,y\in\Z),
\\0\ (\mo\ p^2)&\t{if}\ p\eq2\ (\mo\ 3).\endcases
\endaligned\tag1.10$$
\endproclaim
\Proof. By [M05] and [Su4], we have
$$\sum_{k=0}^{p-1}\f{\bi{4k}{k,k,k,k}}{256^k}\eq c(p)\pmod{p^2}$$
as conjectured in [RV]. (Here we only need the mod $p$ version which
was proved in [M05].)  So (1.9) follows from (1.4).  The author
[Su2] conjectured that
$$\align&\sum_{k=0}^{p-1}\f{\bi{2k}k^3}{16^k}\eq\l(\f{-1}p\r)\sum_{k=0}^{p-1}\f{\bi{2k}k^3}{256^k}
\\\eq&\cases4x^2-2p\ (\mo\ p^2)&\t{if}\ p\eq1\ (\mo\ 3)\ \t{and}\ p=x^2+3y^2\ (x,y\in\Z),
\\0\ (\mo\ p^2)&\t{if}\ p\eq2\ (\mo\ 3).\endcases
\endalign$$
This was confirmed by Z. H. Sun [S3] in the case $p\eq2\pmod3$, and the mod $p$ version in the case $p\eq1\pmod3$
follows from (4)-(5) in Ahlgren [A, Theorem 5]. So we get (1.10) by applying (1.3) with $x=1,1/16$.
\qed

\Remark\ 1.2. The author conjectured that (1.9) also holds modulo $p^2$, and that (1.10)
is also valid modulo $p^2$ in the case $p\eq1\pmod3$.

\proclaim{Corollary 1.3} For any odd prime $p$ and integer $x$, we have
$$\sum_{k=0}^{p-1}(2k+1)A_k(x)\eq p\l(\f xp\r)\ (\mo\ p^2).\tag1.11$$
\endproclaim
\Proof. This follows from (1.5) in the case $n=p$, for,
$p\mid\bi{p+k}{2k+1}$ for every $k=0,\ldots,(p-3)/2$, and $p\mid\bi{2k}k$ for all $k=(p+1)/2,\ldots,p-1$.
\qed

\medskip

We deduce Theorem 1.1(i) from our following result which has its own interest.

\proclaim{Theorem 1.2} Let $p$ be an odd prime and let $x$ be any $p$-adic integer.

{\rm (i)} If $x\eq2k\pmod{p}$ with $k\in\{0,\ldots,(p-1)/2\}$,
then we have
$$\sum_{r=0}^{p-1}(-1)^r\bi xr^2\eq(-1)^k\bi xk\pmod{p^2}.\tag1.12$$

{\rm (ii)} If $x\eq k\pmod{p}$ with $k\in\{0,\ldots,p-1\}$, then
$$\sum_{r=0}^{p-1}\bi xr^2\eq\bi{2x}k\pmod{p^2}.\tag1.13$$
\endproclaim
\Remark\ 1.3. In contrast with (1.12) and (1.13), we recall the following identities (cf. [G, (3.32) and (3.66)]):
$$\sum_{k=0}^{2n}(-1)^k\bi {2n}k^2=(-1)^n\bi {2n}n\quad\t{and}\quad\sum_{k=0}^n\bi nk^2=\bi{2n}n.$$

\proclaim{Corollary 1.4} Let $p$ be an odd prime.

{\rm (i) (Conjectured in [RV] and proved in [M03])} We have
$$\sum_{k=0}^{p-1}\f{\bi{2k}k^2}{16^k}\eq\l(\f{-1}p\r)\pmod{p^2}.$$

{\rm (ii) (Conjectured by the author [Su1] and confirmed in [S2])} If
$p\eq1\pmod4$ and $p=x^2+y^2$ with $x\eq1\pmod4$ and $y\eq0\pmod2$, then
$$\sum_{k=0}^{p-1}\f{\bi{2k}k^2}{(-16)^k}\eq(-1)^{(p-1)/4}\l(2x-\f p{2x}\r)\pmod{p^2}.\tag1.14$$
\endproclaim
\Proof. Since $\bi{-1/2}{r}=\bi{2r}r/(-4)^r$ for all $r=0,1,\ldots$, applying (1.13) with $x=-1/2$ and $k=(p-1)/2$
we immediately get the congruence in part (i).

When $p=x^2+y^2$ with $x\eq1\pmod4$ and $y\eq0\pmod2$, by (1.12) with $x=-1/2$ and
$k=(p-1)/4$ we have
$$\align\sum_{r=0}^{p-1}\f{\bi{2r}r^2}{(-16)^r}\eq&(-1)^{(p-1)/4}\bi{-1/2}{(p-1)/4}=\f{\bi{(p-1)/2}{(p-1)/4}}{4^{(p-1)/4}}
=\f{\bi{(p-1)/2}{(p-1)/4}}{2^{(p-1)/2}}
\\\eq&\f{2^{p-1}+1}{2\times2^{(p-1)/2}}\l(2x-\f p{2x}\r)\pmod{p^2}\ \, \t{(by [CDE] or [BEW, (9.0.2)])}
\\\eq&(-1)^{(p-1)/4}\l(2x-\f p{2x}\r)\pmod{p^2}
\endalign$$
since $((-1)^{(p-1)/4}2^{(p-1)/2}-1)^2\eq 0\ (\mo\ p^2)$.
This proves (1.14). \qed
\proclaim{Corollary 1.5} Let $a_n:=\sum_{k=0}^n\bi nk^2C_k$ for $n=0,1,2,\ldots$, where
$C_k$ denotes the Catalan number $\bi{2k}k/(k+1)=\bi{2k}k-\bi{2k}{k+1}$. Then,
for any odd prime $p$ we have
$$a_1+\cdots+a_{p-1}\eq0\pmod{p^2}.\tag1.15$$
\endproclaim
\Remark\ 1.4. We find no prime $p\ls 5,000$ with $\sum_{k=1}^{p-1}a_k\eq0\pmod{p^3}$ and no composite number $n\ls70,000$
satisfying $\sum_{k=1}^{n-1}a_k\eq0\pmod{n^2}$. We conjecture that (1.15) holds for no composite $p>1$.
\medskip

The author [Su1, Remark 1.2] conjectured that for any prime $p>5$ with $p\eq1\pmod4$ we have
$$\sum_{k=0}^{p^a-1}\f{k^3\bi{2k}k^3}{64^k}\eq0\pmod{p^{2a}}\quad\t{for}\ a=1,2,3,\ldots.$$
This was recently confirmed by Z. H. Sun [S3] in the case $a=1$. Note that
$$\f{k^3\bi{2k}k^3}{64^k}=(-1)^kk^3\bi{-1/2}k^3=\f{(-1)^{k-1}}8\bi{-3/2}{k-1}^3\quad\t{for all}\ k=1,2,3,\ldots.$$
So, for any prime $p>5$ with $p\eq1\pmod4$ we have
$$\sum_{r=0}^{p-1}(-1)^r\bi{-3/2}r^3\eq0\pmod{p^2}.$$
Since $-3/2\eq-2(p+3)/4\pmod{p}$, the result just corresponds to the case $x=-3/2$ of our following general theorem.

\proclaim{Theorem 1.3} Let $p>3$ be a prime and let $x$ be a $p$-adic integer with $x\eq-2k\pmod{p}$ for some
$k\in\{1,\ldots,\lfloor(p-1)/3\rfloor\}$. Then we have
$$\sum_{r=0}^{p-1}(-1)^r\bi xr^3\eq0\pmod{p^2}.\tag1.16$$
\endproclaim

Similar to Ap\'ery numbers, the central Delannoy numbers (see [CHV]) are defined by
$$D_n=\sum_{k=0}^n\bi {n+k}{2k}\bi{2k}k=\sum_{k=0}^n\bi nk\bi{n+k}k\ (n\in\N).$$
Such numbers arise naturally in many enumeration problems in combinatorics (cf. Sloane [S]); for example, $D_n$
is the number of lattice paths from $(0,0)$ to $(n,n)$ with steps $(1,0),(0,1)$ and $(1,1)$.

Now we give our result on central Delannoy numbers.

\proclaim{Theorem 1.4} Let $p>3$ be a prime. Then
$$\sum_{k=0}^{p-1}D_k\eq\l(\f{-1}p\r)-p^2E_{p-3}\ (\mo\ p^3),\tag1.17$$
We also have
$$\sum_{k=0}^{p-1}(2k+1)(-1)^kD_k\eq p-\f 7{12}p^4B_{p-3}\ (\mo\ p^5)\tag1.18$$
and
$$\sum_{k=0}^{p-1}(2k+1)D_k\eq p+2p^2q_p(2)-p^3q_p(2)^2\ (\mo\ p^4),\tag1.19$$
where $q_p(2)$ denotes the Fermat quotient $(2^{p-1}-1)/p$.
\endproclaim

\Remark\ 1.5. In [Su3] the author determined $\sum_{k=1}^{p-1}D_k/k$ and $\sum_{k=1}^{p-1}D_k/k^2$
modulo an odd prime $p$.
\medskip

In the next section we will show Theorems 1.1-1.2 and Corollary 1.5.
Section 3 is devoted to our proofs of Theorems 1.3 and 1.4.
In Section 4 we are going to raise some related conjectures for further research.

\heading{2. Proofs of Theorems 1.1-1.2 and Corollary 1.5}\endheading

We first prove Theorem 1.2.

\medskip
\noindent{\it Proof of Theorem 1.2}. (i) We now consider the first part of Theorem 1.2. Set
$$f_k(y):=\sum_{r=0}^{p-1}(-1)^r\bi {2k+py}r^2\ \ \quad\t{for}\ k\in\N.\tag2.1$$
We want to prove that
$$f_k(y)\eq(-1)^k\bi{2k+py}k\pmod {p^2}\tag2.2$$
for any $p$-adic integer $y$ and $k\in\{0,1,\ldots,(p-1)/2\}$.

Applying the Zeilberger algorithm (cf. [PWZ]) via {\tt Mathematica 7}, we find that
$$\aligned&(py+2k+2)f_{k+1}(y)+4(py+2k+1)f_k(y)
\\=&\f{(p(y-1)+2k+3)^2 F_k(y)}{(py+2k+1)(py+2k+2)^2}\bi{py+2k+2}{p-1}^2,
\endaligned\tag2.3$$
where
$$F_k(y)=14+34k+20k^2-10p-12kp+2p^2+17py+20kpy-6p^2y+5p^2y^2.$$

Now fix a $p$-adic integer $y$. Observe that
$$\align f_{(p-1)/2}(y)=&\sum_{r=0}^{p-1}(-1)^r\bi{p-1+py}r^2=\sum_{r=0}^{p-1}(-1)^r\prod_{0<s\ls r}\l(1-\f{p(y+1)}s\r)^2
\\\eq&\sum_{r=0}^{p-1}(-1)^r\(1-\sum_{0<s\ls r}\f{2p(y+1)}s\)=1-\sum_{r=1}^{p-1}(-1)^r\sum_{s=1}^r\f{2p(y+1)}s
\\=&1-2p(y+1)\sum_{s=1}^{p-1}\f1s\sum_{r=s}^{p-1}(-1)^r=1-p(y+1)\sum_{j=1}^{(p-1)/2}\f1j
\\\eq&(-1)^{(p-1)/2}\bi{p-1+py}{(p-1)/2}\pmod{p^2}.
\endalign$$
For each $k\in\{0,\ldots,(p-3)/2\}$, clearly
$py+2k+1,py+2k+2\not\eq0\pmod{p}$,  and also
$$(p(y-1)+2k+3)^2\bi{py+2k+2}{p-1}^2\eq0\pmod{p^2}$$
since $\bi{py+2k+2}{p-1}=\f p{py+2k+3}\bi{py+2k+3}p\eq0\pmod p$ if $0\ls k<(p-3)/2$.
Thus, by (2.3) we have
$$f_k(y)\eq-\f{py+2k+2}{4(py+2k+1)}f_{k+1}(y)\pmod{p^2}\quad\t{for}\ k=0,\ldots,\f{p-3}2.$$
If $0\ls k<(p-1)/2$ and
$$f_{k+1}(y)\eq(-1)^{k+1}\bi{2(k+1)+py}{k+1}\pmod{p^2},$$
then
$$\align f_k(y)\eq &-\f{py+2k+2}{4(py+2k+1)}(-1)^{k+1}\bi{2(k+1)+py}{k+1}
\\=&\f{(-10^k(py+2k+2)^2}{4(k+1)(py+k+1)}\bi{2k+py}k\eq(-1)^k\bi{2k+py}k\pmod{p^2}.
\endalign$$
Therefore (2.2) holds for all $k=0,1,\ldots,(p-1)/2$. This proves Theorem 1.2(i).

(ii) The second part of Theorem 1.2 can be proved in a similar way. Here we mention that
if we define
$$g_k(y):=\sum_{r=0}^{p-1}\bi{k+py}r^2\ \ \t{for}\ k\in\N\tag2.4$$
then by the Zeilberger algorithm (cf. [PWZ]) we have the recursion
$$\align&(py+k+1)g_{k+1}(y)-2(2py+2k+1)g_k(y)
\\=&-\f{(p(y-1)+k+2)^2(3py-2p+3k+3)}{(py+k+1)^2}\bi{py+k+1}{p-1}^2.
\endalign$$
It follows that if $k\in\{0,\ldots,p-2\}$ and $y$ is a $p$-adic integer then
$$\aligned &g_{k+1}(y)\eq\bi{2(k+1)+2py}{k+1}\pmod{p^2}
\\\Longrightarrow\ &g_k(y)\eq\bi{2k+2py}{k}\pmod{p^2}.
\endaligned\tag2.5$$
In view of this, we have the second part of Theorem 1.2 by induction.

The proof of Theorem 1.2 is now complete. \qed

\medskip\noindent{\it Proof of Corollary 1.5}. Observe that
$$\sum_{n=0}^{p-1}a_n=\sum_{k=0}^{p-1}C_k\sum_{n=k}^{p-1}\bi nk^2=\sum_{k=0}^{p-1}C_k\sum_{j=0}^{p-1-k}\bi{k+j}k^2.$$
If $0\ls k\ls p-1$ and $p-k\ls j\ls p-1$, then
$$\bi{k+j}k=\f{(k+j)!}{k!j!}\eq0\pmod p.$$
Therefore
$$\sum_{n=0}^{p-1}a_n\eq\sum_{k=0}^{p-1}C_k\sum_{j=0}^{p-1}\bi{k+j}k^2=\sum_{k=0}^{p-1}C_k\sum_{j=0}^{p-1}\bi{x_k}j^2,$$
where $x_k=-k-1\eq p-1-k\pmod{p}$. Applying Theorem 1.2(ii) we get
$$\sum_{n=0}^{p-1}a_n\eq\sum_{k=0}^{p-1}C_k\bi{2x_k}{p-1-k}=\sum_{k=0}^{p-1}(-1)^k\bi{p+k}{2k+1}C_k\pmod{p^2}.$$
So it suffices to show that for any $n\in\Z^+$ we have
$$\sum_{k=0}^{n-1}(-1)^k\bi{n+k}{2k+1}C_k=1.\tag2.6$$

We prove (2.6) by induction. Clearly, (2.6) holds for $n=1$.
Let $n$ be any positive integer. By the Chu-Vandermonde identity
$$\sum_{k=0}^n\bi xk\bi y{n-k}=\bi{x+y}n$$
(see, e.g., [GKP, p.\,169]), we have
$$\sum_{k=0}^{n-1}\bi{n+1}{k+1}\bi{n+k}k(-1)^k=\sum_{k=0}^{n-1}\bi{n+1}{n-k}\bi{-n-1}k=-\bi{-n-1}n.$$
Thus
$$\align &\sum_{k=0}^n(-1)^k\bi{n+1+k}{2k+1}C_k-\sum_{k=0}^{n-1}(-1)^k\bi{n+k}{2k+1}C_k
\\=&(-1)^nC_n+\sum_{k=0}^{n-1}(-1)^k\bi{n+k}{2k}C_k
\\=&(-1)^nC_n+\f1{n+1}\sum_{k=0}^{n-1}\bi{n+1}{k+1}\bi{n+k}k(-1)^k
\\=&(-1)^nC_n-\f1{n+1}\bi{-n-1}n=0.
\endalign$$
This concludes the induction step. We are done. \qed

\medskip
Now we can apply Theorem 1.2 to deduce the first part of Theorem 1.1.

\medskip\noindent{\it Proof of Theorem} 1.1(i).  Let $\ve\in\{\pm1\}$.
Then
$$\align\sum_{m=0}^{p-1}\ve^mA_m(x)=&\sum_{m=0}^{p-1}\ve^m\sum_{k=0}^{m}\bi{m+k}{2k}^2\bi{2k}k^2x^k
\\=&\sum_{k=0}^{p-1}\bi{2k}k^2x^k\sum_{m=k}^{p-1}\ve^m\bi{m+k}{2k}^2
\\=&\sum_{k=0}^{p-1}\bi{2k}k^2x^k\sum_{r=0}^{p-1-k}\ve^{k+r}\bi{2k+r}{r}^2
\\=&\sum_{k=0}^{p-1}\bi{2k}k^2\ve^kx^k\sum_{r=0}^{p-1-k}\ve^{r}\bi{p-1-2k-p}{r}^2
\endalign$$
Set $n=(p-1)/2$. Clearly $\bi{2k}k\eq0\pmod p$ for $k=n+1,\ldots,p-1$, and
$$\bi{p-1-2k-p}r\eq\bi{p-1-2k}r=0\pmod{p}$$
if $0\ls k\ls n$ and $p-1-2k<r\ls p-1$. Therefore
$$\sum_{m=0}^{p-1}\ve^mA_m(x)\eq\sum_{k=0}^{n}\bi{2k}k^2\ve^kx^k\sum_{r=0}^{p-1}\ve^r\bi{2(n-k)-p}r^2\pmod{p^2}.$$
Similarly,
$$\align\sum_{m=0}^{p-1}\ve^mW_m(\ve x)=&\sum_{m=0}^{p-1}\ve^m\sum_{k=0}^{\lfloor m/2\rfloor}\bi{m}{2k}^2\bi{2k}k^2(\ve x)^k
\\=&\sum_{k=0}^{n}\bi{2k}k^2\ve^k x^k\sum_{m=2k}^{p-1}\ve^m\bi{m}{2k}^2
\\=&\sum_{k=0}^{n}\bi{2k}k^2\ve^k x^k\sum_{r=0}^{p-1-2k}\ve^{2k+r}\bi{2k+r}{r}^2
\\\eq&\sum_{k=0}^{n}\bi{2k}k^2\ve^kx^k\sum_{r=0}^{p-1}\ve^{r}\bi{2(n-k)-p}{r}^2\pmod{p^2}.
\endalign$$
So we have
$$\sum_{m=0}^{p-1}\ve^mA_m(x)\eq\sum_{m=0}^{p-1}\ve^mW_m(\ve x)
\eq\sum_{k=0}^{n}\bi{2k}k^2\ve^kx^k S_k(\ve)\pmod{p^2},\tag2.7$$
where
$$S_k(\ve):=\sum_{r=0}^{p-1}\ve^{r}\bi{2(n-k)-p}{r}^2.$$

 Applying Theorem 1.2(i) we get
$$\align S_k(-1)\eq&(-1)^{n-k}\bi{2(n-k)-p}{n-k}=(-1)^{n-k}\bi{-2k-1}{n-k}
\\=&\bi{n+k}{n-k}=\bi{n+k}{2k}\eq\f{\bi{2k}k}{(-16)^k}\pmod{p^2}.
\endalign$$
(The last congruence can be easily deduced, see. e.g., [S2, Lemma 2.2].)
Combining this with (2.7) in the case $\ve=-1$ we immediately obtain (1.3).

In view of Theorem 1.2(ii),
$$S_k(1)\eq\bi{4(n-k)-2p}{2(n-k)}\pmod{p^2}.$$
Recall that $\bi{n+k}{n-k}(-16)^k\eq\bi{2k}k\pmod{p^2}.$
So, in view of (2.7) with $\ve=1$, we have
$$\align\sum_{m=0}^{p-1}A_m(x)\eq&\sum_{m=0}^{p-1}W_m(x)
\eq\sum_{k=0}^{n}\bi{n+k}{n-k}^2(-16)^{2k}x^k\bi{4(n-k)-2p}{2(n-k)}
\\=&\sum_{j=0}^n\bi{n+(n-j)}j^2256^{n-j}x^{n-j}\bi{4j-2p}{2j}
\\=&16^{p-1}\sum_{k=0}^n\f{\bi{4k-2p}{2k}\bi{2k-p}k^2}{256^k}x^{n-k} \pmod{p^2}
\endalign$$
If $x$ is a $p$-adic integer with $x\not\eq0\pmod p$, then
$$\align &16^{p-1}\sum_{k=0}^n\f{\bi{4k-2p}{2k}\bi{2k-p}k^2}{256^k}x^{n-k}
\\\eq&\l(\f xp\r)\sum_{k=0}^n\f{\bi{4k}{2k}\bi{2k}k^2}{(256x)^k}
\eq\l(\f xp\r)\sum_{k=0}^{p-1}\f{\bi{4k}{k,k,k,k}}{(256x)^k}\pmod p,
\endalign$$
and therefore (1.4) holds. \qed

\proclaim{Lemma 2.1} Let $k\in\N$. Then, for any $n\in\Z^+$ we have the identity
$$\sum_{m=0}^{n-1}(2m+1)\bi{m+k}{2k}^2=\f{(n-k)^2}{2k+1}\bi{n+k}{2k}^2.\tag 2.8$$
\endproclaim
\Proof. Obviously (2.8) holds when $n=1$.

 Now assume that $n>1$ and (2.8) holds. Then
$$\align&\sum_{m=0}^{n}(2m+1)\bi{m+k}{2k}^2
\\=&\f{(n-k)^2}{2k+1}\bi{n+k}{2k}^2+(2n+1)\bi{n+k}{2k}^2
\\=&\f{(n+k+1)^2}{2k+1}\bi{n+k}{2k}^2=\f{(n+1-k)^2}{2k+1}\bi{(n+1)+k}{2k}^2.
\endalign$$

Combining the above, we have proved the desired result by induction. \qed

\proclaim{Lemma 2.2} Let $p>3$ be a prime. Then
$$\sum^{p-1}\Sb k=0\\k\not=(p-1)/2\endSb\f{(-1)^k}{2k+1}\eq-pE_{p-3}\ (\mo\ p^2).\tag2.9$$
\endproclaim
\Proof. Observe that
$$\align \sum^{p-1}\Sb k=0\\k\not=(p-1)/2\endSb\f{(-1)^k}{2k+1}
=&\f12\sum^{p-1}\Sb k=0\\k\not=(p-1)/2\endSb\(\f{(-1)^k}{2k+1}+\f{(-1)^{p-1-k}}{(2(p-1-k)+1)}\)
\\=&-p\sum^{p-1}\Sb k=0\\k\not=(p-1)/2\endSb\f{(-1)^k}{(2k+1)(2k+1-2p)}
\\\eq&-\f p4\sum_{k=0}^{p-1}(-1)^k\l(k+\f12\r)^{p-3}\ (\mo\ p^2).
\endalign$$
So we have reduced (2.9) to the following congruence
$$\sum_{k=0}^{p-1}(-1)^k\l(k+\f12\r)^{p-3}\eq 4E_{p-3}\ (\mo\ p).\tag2.10$$

Recall that the Euler polynomial of degree $n$ is defined by
$$E_n(x)=\sum_{k=0}^n\bi nk\f{E_k}{2^k}\l(x-\f12\r)^{n-k}.$$
It is well known that
$$E_n(x)+E_n(x+1)=2x^n.$$
Thus
$$\align &2\sum_{k=0}^{p-1}(-1)^k\l(k+\f12\r)^{p-3}
\\=&\sum_{k=0}^{p-1}\l((-1)^kE_{p-3}\l(k+\f12\r)-(-1)^{k+1}E_{p-3}\l(k+1+\f12\r)\r)
\\=&E_{p-3}\l(\f12\r)-(-1)^pE_{p-3}\l(p+\f12\r)
\\\eq&2E_{p-3}\l(\f12\r)=2\f{E_{p-3}}{2^{p-3}}\eq8E_{p-3}\ (\mo\ p)
\endalign$$
and hence (2.10) follows. We are done. \qed

\medskip
For each $m=1,2,3,\ldots$ those  rational numbers
$$H_n^{(m)}:=\sum_{0<k\ls n}\f1{k^m}\ \ (n=0,1,2,\ldots)$$
are called harmonic numbers of order $m$. We simply write $H_n$ for $H_n^{(1)}$.
A well-known theorem of Wolstenholme asserts that $H_{p-1}\eq0\pmod{p^2}$ and $H_{p-1}^{(2)}\eq0\pmod p$ for any prime $p>3$.

\proclaim{Lemma 2.3} Let $p>3$ be a prime. Then
$$\sum_{k=0}^{(p-3)/2}\f{H_k^{(2)}}{2k+1}\eq-\f 74B_{p-3}\pmod p.\tag2.11$$
\endproclaim
\Proof. Clearly,
$$\sum_{k=1}^{p-1}\f1{k^3}=\sum_{k=1}^{(p-1)/2}\l(\f1{k^3}+\f1{(p-k)^3}\r)\eq0\pmod{p}.$$
By [ST, (5.4)], $\sum_{k=1}^{p-1}H_k/k^2\eq B_{p-3}\pmod{p}$. Therefore
$$\align\sum_{k=1}^{p-1}\f{H_k^{(2)}}k=&\sum_{k=1}^{p-1}\f1k\sum_{j=1}^k\f1{j^2}=\sum_{j=1}^{p-1}\f{H_{p-1}-H_{j-1}}{j^2}
\\\eq&-\sum_{k=1}^{p-1}\f{H_k}{k^2}+\sum_{k=1}^{p-1}\f1{k^3}\eq-B_{p-3}\pmod p.
\endalign$$
On the other hand,
$$\align&\sum_{k=1}^{p-1}\f{H_k^{(2)}}k=\sum_{k=1}^{(p-1)/2}\l(\f{H_k^{(2)}}k+\f{H_{p-k}^{(2)}}{p-k}\r)
\\\eq&\sum_{k=1}^{(p-1)/2}\l(\f{H_k^{(2)}}k+\f{H_{p-1}^{(2)}-H_{k-1}^{(2)}}{-k}\r)\eq2\sum_{k=1}^{(p-1)/2}\f{H_k^{(2)}}k-H_{(p-1)/2}^{(3)}
\pmod p.\endalign$$
It is known (see, e.g., [S1, Corollary 5.2]) that
$$H_{(p-1)/2}^{(3)}=\sum_{k=1}^{(p-1)/2}\f1{k^3}\eq-2B_{p-3}\pmod p.$$
So we have
$$\sum_{k=1}^{(p-1)/2}\f{H_k^{(2)}}k\eq\f12\(\sum_{k=1}^{p-1}\f{H_k^{(2)}}k+H_{(p-1)/2}^{(3)}\)
\eq\f{-B_{p-3}-2B_{p-3}}2=-\f 32B_{p-3}\pmod p.$$

Clearly
$$H_{(p-1)/2}^{(2)}\eq\f12\sum_{k=1}^{(p-1)/2}\l(\f1{k^2}+\f1{(p-k)^2}\r)=\f12H_{p-1}^{(2)}\eq0\pmod{p}.$$
Observe that
$$\align \sum_{k=0}^{(p-3)/2}\f{H_k^{(2)}}{2k+1}\eq&-\sum_{k=0}^{(p-3)/2}\f{H_k^{(2)}}{p-1-2k}=-\sum_{k=1}^{(p-1)/2}\f{H_{(p-1)/2-k}^{(2)}}{2k}
\\\eq&-\f12\sum_{k=1}^{(p-1)/2}\f1k\(H_{(p-1)/2}^{(2)}-\sum_{j=0}^{k-1}\f1{((p-1)/2-j)^2}\)
\\\eq&2\sum_{k=1}^{(p-1)/2}\f1k\sum_{j=0}^{k-1}\f1{(2j+1)^2}\eq2\sum_{k=1}^{(p-1)/2}\f1k\(H_{2k}^{(2)}-\sum_{j=1}^k\f1{(2j)^2}\)
\\=&4\sum_{k=1}^{(p-1)/2}\f{H_{2k}^{(2)}}{2k}-\f12\sum_{k=1}^{(p-1)/2}\f{H_k^{(2)}}k\pmod p
\endalign$$
and
$$\align\sum_{k=1}^{(p-1)/2}\f{H_{2k}^{(2)}}{2k}=&\sum^{p-1}\Sb k=1\\2\mid k\endSb\f1k\sum_{j=1}^k\f1{j^2}
=\sum^{p-1}\Sb k=1\\2\mid k\endSb\f1{k^3}+\sum\Sb 1\ls j<k\ls p-1\\2\mid k\endSb\f1{j^2k}
\\\eq&\f18H_{(p-1)/2}^{(3)}-\f 38B_{p-3}\quad\t{(by Pan [P, (2.4)])}
\\\eq&\f18(-2B_{p-3})-\f 38B_{p-3}=-\f 58B_{p-3}\pmod p.
\endalign$$
So we finally get
$$\sum_{k=0}^{(p-3)/2}\f{H_k^{(2)}}{2k+1}\eq4\l(-\f 58B_{p-3}\r)-\f12\l(-\f 32B_{p-3}\r)=-\f 74B_{p-3}\pmod p.$$
This concludes the proof of (2.11). \qed

\medskip\noindent{\it Proof of Theorem} 1.1(ii). (i) Let $n$ be any positive integer. Then
$$\align \sum_{m=0}^{n-1}(2m+1)A_m(x)=&\sum_{m=0}^{n-1}(2m+1)\sum_{k=0}^m\bi{m+k}{2k}^2\bi{2k}k^2x^k
\\=&\sum_{k=0}^{n-1}\bi{2k}k^2x^k\sum_{m=0}^{n-1}(2m+1)\bi{m+k}{2k}^2
\\=&\sum_{k=0}^{n-1}\bi{2k}k^2x^k\f{(n-k)^2}{2k+1}\bi{n+k}{2k}^2\ \ (\t{by (2.8)})
\\=&\sum_{k=0}^{n-1}\f{(n-k)^2}{2k+1}\bi n k^2\bi{n+k}k^2x^k.
\endalign$$
Since
$$(n-k)\bi nk=n\bi{n-1}k\ \quad\t{for all}\ \ k=0,\ldots,n-1,$$
we have
$$\align\f1n\sum_{m=0}^{n-1}(2m+1)A_m(x)=&\sum_{k=0}^{n-1}\bi{n-1}k\f {n-k}{2k+1}\bi nk\bi{n+k}k^2x^k
\\=&\sum_{k=0}^{n-1}\bi{n-1}k\f {n-k}{2k+1}\bi{n+k}{2k}\bi {2k}k\bi{n+k}kx^k
\\=&\sum_{k=0}^{n-1}\bi{n-1}k\bi{n+k}k\bi{n+k}{2k+1}\bi{2k}kx^k.
\endalign$$
This proves (1.5).

Now we fix a prime $p>3$. By the above,
$$\sum_{m=0}^{p-1}(2m+1)A_m(x)=\sum_{k=0}^{p-1}\f{p^2}{2k+1}\bi{p-1}k^2\bi{p+k}k^2x^k.\tag2.12$$
For $k\in\{0,\ldots,p-1\}$, clearly
$$\align\bi{p-1}k^2\bi{p+k}k^2=&\prod_{0<j\ls k}\l(\f{p-j}j\cdot\f{p+j}j\r)^2=\prod_{0<j\ls k}\l(1-\f{p^2}{j^2}\r)^2
\\\eq&\prod_{0<j\ls k}\l(1-\f{2p^2}{j^2}\r)\eq1-2p^2H_k^{(2)}\pmod{p^4}.
\endalign$$
Thus (2.12) implies that
$$\sum_{m=0}^{p-1}(2m+1)A_m(x)=\sum_{k=0}^{p-1}\f{p^2}{2k+1}\l(1-2p^2H_k^{(2)}\r)x^k\pmod{p^5}.\tag2.13$$

Since $H_{(p-1)/2}^{(2)}\eq0\pmod p$, taking $x=-1$ in (2.13) and applying (2.9) we obtain
$$\sum_{m=0}^{p-1}(2m+1)A_m(-1)\eq\sum_{k=0}^{p-1}\f{p^2(-1)^k}{2k+1}
\eq \f{p^2(-1)^{(p-1)/2}}{2(p-1)/2+1}-p^3E_{p-3}\pmod{p^4}$$
and hence (1.7) holds.

Now we prove (1.6). In view of (2.13) with $x=1$, we have
$$\align\sum_{m=0}^{p-1}(2m+1)A_m\eq&\f{p^2}{2(p-1)/2+1}\l(1-2p^2H_{(p-1)/2}^{(2)}\r)
\\&+p^2\sum_{k=0}^{(p-3)/2}\l(\f{1-2p^2H_k^{(2)}}{2k+1}+\f{1-2p^2H_{p-1-k}^{(2)}}{2(p-1-k)+1}\r)
\\=&p-2p^3H_{(p-1)/2}^{(2)}+2p^3\sum_{k=0}^{(p-3)/2}\f{2p+2k+1}{(2k+1)(4p^2-(2k+1)^2)}
\\&-2p^4\sum_{k=0}^{(p-3)/2}\(\f{H_k^{(2)}}{2k+1}+\f{H_{p-1}^{(2)}-\sum_{0<j\ls k}(p-j)^{-2}}{2p-(2k+1)}\)
\\\eq&p-2p^3H_{(p-1)/2}^{(2)}-4p^4\sum_{k=0}^{(p-3)/2}\f1{(2k+1)^3}
\\&-2p^3\sum_{k=0}^{(p-3)/2}\f1{(2k+1)^2}-4p^4\sum_{k=0}^{(p-3)/2}\f{H_k^{(2)}}{2k+1}\pmod{p^5}.
\endalign$$
By [S1, Corollaries 5.1 and 5.2],
$$H_{p-1}^{(2)}\eq\f23pB_{p-3}\pmod{p^2},\quad H_{(p-1)/2}^{(2)}\eq\f 73pB_{p-3}\pmod{p^2},$$
$$\sum_{k=0}^{(p-3)/2}\f1{(2k+1)^2}=H_{p-1}^{(2)}-\f{H_{(p-1)/2}^{(2)}}4\eq\f p{12}B_{p-3}\pmod{p^2},$$
and
$$\sum_{k=0}^{(p-3)/2}\f1{(2k+1)^3}=H_{p-1}^{(3)}-\f{H_{(p-1)/2}^{(3)}}8\eq0-\f{-2B_{p-3}}8=\f{B_{p-3}}4\pmod p.$$
Combining these with Lemma 2.3, we finally obtain
$$\align\sum_{k=0}^{p-1}(2k+1)A_k\eq &p-2p^3\f 73pB_{p-3}-4p^4\f{B_{p-3}}4-2p^3\f p{12}B_{p-3}-4p^4\l(-\f 74B_{p-3}\r)
\\=&p+\f 76p^4B_{p-3}\pmod{p^5}
\endalign$$

So far we have proved the second part of Theorem 1.1. \qed

\medskip

Part (iii) of Theorem 1.1 is easy.

\medskip
\noindent{\it Proof of Theorem} 1.1(iii). As $A_0=1$ and $A_1=3$, the desired congruence with $p=2$ holds trivially.

Below we assume that $p$ is an odd prime. If $k\in\{0,1\ldots,p-1\}$, then
$$\align&A_{p-1-k}=\sum_{j=0}^{p-1}\bi{(p-1-k)+j}{2j}^2\bi{2j}j^2
\\\eq&\sum_{j=0}^{p-1}\bi{j-k-1}{2j}^2\bi{2j}j^2=\sum_{j=0}^k\bi{j+k}{2j}^2\bi{2j}j^2=A_k\ (\mo\ p)
\endalign$$
Thus
$$\align\sum_{k=0}^{p-1}(2k+1)\ve^kA_k^m&=\sum_{k=0}^{p-1}(2(p-1-k)+1)\ve^{p-1-k}A_{p-1-k}^m
\\&\eq-\sum_{k=0}^{p-1}(2k+1)\ve^kA_k^m\ (\mo\ p)
\endalign$$
and hence we have the desired congruence. \qed

\heading{3. Proofs of Theorems 1.3 and 1.4}\endheading

\medskip
\noindent{\it Proof of Theorem 1.3}. Define
$$w_k(y):=\sum_{r=0}^{p-1}(-1)^r\bi{py-2k}r^3\quad\t{for}\ k\in\N.\tag3.1$$
We want to show that $w_k(y)\eq0\pmod{p^2}$ for any $p$-adic integer $y$ and
$k\in\{1,\ldots,\lfloor(p-1)/3\rfloor\}$.

By the Zeilberger algorithm (cf. [PWZ]), for $k=0,1,2,\ldots$ we have
$$\aligned&(py-2k)^2 w_k(y)+3(3py-2(3k+1))(3py-2(3k+2))w_{k+1}(y)
\\&\qquad=\f{P(k,p,y)(p(1-y)+2k-1)^3}{(py-2k)^3(py-2k-1)^3}\bi{py-2k}{p-1}^3
\endaligned\tag3.2$$
where $P(k,p,y)$ is a suitable polynomial in $k,p,y$ with integer coefficients such that
$P(0,p,y)\eq0\pmod{p^2}$. (Here we omit the explicit expression of $P(k,p,y)$ since it is complicated.)
Note also that
$$w_1(0)=\sum_{r=0}^{p-1}(-1)^r\bi{-2}r^3=\sum_{r=0}^{p-1}(r+1)^3=\f{p^2(p+1)^2}4\eq0\pmod{p^2}.$$

 Fix a $p$-adic integer $y$. If $y\not=0$, then (3.2) with $k=0$ yields
$$\align&3(3py-2)(3py-4)w_1(y)
\\\eq&\f{P(0,p,y)(p(1-y)-1)^3}{(py)^3(py-1)^3}\l(\f{py}{p-1}\bi{p(y-1)+p-1}{p-2}\r)^3
\eq0\pmod{p^2}
\endalign$$
and hence $w_1(y)\eq0\pmod{p^2}$. If $1<k+1\ls\lfloor(p-1)/3\rfloor$, then by (3.2) we have
$$(py-2k)^2 w_k(y)+3(3py-2(3k+1))(3py-2(3k+2))w_{k+1}(y)\eq0\pmod{p^3}$$
since
$$\bi{py-2k}{p-1}=\f p{py-2k+1}\bi{py-2k+1}p\eq0\pmod{p}.$$
Thus, when $1<k+1\ls \lfloor(p-1)/3\rfloor$ we have
$$w_k(y)\eq0\pmod{p^2}\ \Longrightarrow\ w_{k+1}(y)\eq0\pmod{p^2}.$$
So, by induction, $w_k(y)\eq0\pmod{p^2}$ for all $k=1,\ldots,\lfloor(p-1)/3\rfloor$.

In view of the above, we have completed the proof of Theorem 1.3. \qed

\proclaim{Lemma 3.1} Let $n\in\N$. Then we have
$$\sum_{k=0}^n\bi{x+k-1}k=\bi{x+n}n.\tag3.3$$
\endproclaim
\Proof. By the Chu-Vandermonde identity (see, e.g., [GKP, p.\,169]),
$$\sum_{k=0}^n\bi {-x}k\bi {-1}{n-k}=\bi{-x-1}n$$
which is equivalent to (3.3). Of course, it is easy to prove (3.3) by induction. \qed

\medskip\noindent{\it Proof of Theorem 1.4}.
 (i) Observe that
$$\align \sum_{n=0}^{p-1}D_n=&\sum_{n=0}^{p-1}\sum_{k=0}^n\bi{n+k}{2k}\bi{2k}k
=\sum_{k=0}^{p-1}\bi{2k}k\sum_{n=k}^{p-1}\bi{n+k}{2k}
\\=&\sum_{k=0}^{p-1}\bi{2k}k\sum_{j=0}^{p-1-k}\bi{j+2k}j
\\=&\sum_{k=0}^{p-1}\bi{2k}k\bi{2k+1+p-1-k}{p-1-k}\ (\t{by Lemma 3.1})
\\=&\sum_{k=0}^{p-1}\bi{2k}k\bi{p+k}{2k+1}=\sum_{k=0}^{p-1}\f{k+1}{2k+1}\bi{2k+1}k\bi{p+k}{2k+1}
\endalign$$ and thus
$$\align\sum_{n=0}^{p-1}D_n
=&\sum_{k=0}^{p-1}\f{k+1}{2k+1}\bi{p+k}k\bi p{k+1}=p+\sum_{k=1}^{p-1}\f p{2k+1}\bi{p-1}k\bi{p+k}k.
\endalign$$
For $k=1,\ldots,p-1$ we clearly have
$$\bi{p-1}k\bi{p+k}k=(-1)^k\prod_{j=1}^k\l(1-\f{p^2}{j^2}\r)\eq(-1)^k(1-p^2H_k^{(2)})\ (\mo\ p^4);\tag3.4$$
in particular,
$$\bi{p-1}{(p-1)/2}\bi{p+(p-1)/2}{(p-1)/2}\eq(-1)^{(p-1)/2}=\l(\f{-1}p\r)\ (\mo\ p^3)$$
since $H_{(p-1)/2}^{(2)}\eq0\pmod p$.
Therefore
$$\align\sum_{n=0}^{p-1}D_n\eq& \sum^{p-1}\Sb k=0\\k\not=(p-1)/2\endSb\f p{2k+1}(-1)^k+\l(\f{-1}p\r)
\\\eq&\l(\f{-1}p\r)-p^2E_{p-3}\ (\mo\ p^3)\ \ (\t{by (2.9)}).
\endalign$$
This proves (1.17).

(ii) Now we prove (1.18) and (1.19).

Let $n$ be any positive integer. Then
$$\align \sum_{m=0}^{n-1}(2m+1)(-1)^mD_m
=&\sum_{m=0}^{n-1}(2m+1)(-1)^m\sum_{k=0}^m\bi{m+k}{2k}\bi{2k}k
\\=&\sum_{k=0}^{n-1}\bi{2k}k\sum_{m=0}^{n-1}(2m+1)(-1)^m\bi{m+k}{2k}
\endalign$$
By induction, we have the identity
$$\sum_{m=0}^{n-1}(2m+1)(-1)^m\bi{m+k}{2k}=(-1)^n(k-n)\bi{n+k}{2k}.\tag3.5$$
Thus
$$\align \sum_{m=0}^{n-1}(2m+1)(-1)^mD_m=&(-1)^{n-1}\sum_{k=0}^{n-1}\bi{2k}k(n-k)\bi{n+k}{2k}
\\=&(-1)^{n-1}\sum_{k=0}^{n-1}(n-k)\bi nk\bi{n+k}k
\\=&(-1)^{n-1}n\sum_{k=0}^{n-1}\bi{n-1}k\bi{n+k}k.
\endalign$$
Similarly,
$$\align \sum_{m=0}^{n-1}(2m+1)D_m
=&\sum_{m=0}^{n-1}(2m+1)\sum_{k=0}^m\bi{m+k}{2k}\bi{2k}k
\\=&\sum_{k=0}^{n-1}\bi{2k}k\sum_{m=0}^{n-1}(2m+1)\bi{m+k}{2k}
\\=&n\sum_{k=0}^{n-1}C_k (n-k)\bi{n+k}{2k}=\sum_{k=0}^{n-1}\f{n^2}{k+1}\bi{n-1}k\bi{n+k}k.
\endalign$$

In view of (3.4) and the above,
$$\align &\f1p\sum_{m=0}^{p-1}(2m+1)(-1)^mD_m=\sum_{k=0}^{p-1}\bi{p-1}k\bi{p+k}k
\\\eq&\sum_{k=0}^{p-1}(-1)^k-p^2\sum_{k=1}^{p-1}\sum_{0<j\ls k}\f{(-1)^k}{j^2}
=1-p^2\sum_{j=1}^{p-1}\f1{j^2}\sum_{k=j}^{p-1}(-1)^k
\\\eq&1-p^2\sum_{i=1}^{(p-1)/2}\f1{(2i)^2}=1-\f{p^2}4H_{(p-1)/2}^{(2)}\eq1-\f 7{12}p^3B_{p-3}\ (\mo\ p^4)
\endalign$$
and hence (1.18) holds. Similarly,
$$\align &\f1p\sum_{m=0}^{p-1}(2m+1)D_m=\sum_{k=0}^{p-1}\f p{k+1}\bi{p-1}k\bi{p+k}k
\\\eq&\bi{p+(p-1)}{p-1}+p\sum_{k=0}^{p-2}\f{(-1)^k}{k+1}\l(1-p^2H_k^{(2)}\r) \ (\mo\ p^5)
\\\eq&\bi{2p-1}{p-1}-p\sum_{k=1}^{p-1}\f{1+(-1)^k}k\eq1-pH_{(p-1)/2}\ (\mo\ p^3).
\endalign$$
(We have employed  Wolstenholme's congruences
$\bi{2p-1}{p-1}\eq1\ (\mo\ p^3)$ and $H_{p-1}\eq0\ (\mo\ p^2)$.)
To obtain (1.19) it suffices to apply Lehmer's congruence (cf. [L])
$$H_{(p-1)/2}\eq-2q_p(2)+p\,q_p^2(2)\pmod{p^2}.$$

 The proof of Theorem 1.4 is now complete. \qed

\heading{4. Some related conjectures}\endheading

Our following conjecture was motivated by Theorem 1.1(i).

\proclaim{Conjecture 4.1} Let $p>3$ be a prime.

{\rm (i)} If $p\eq1\pmod3$, then
$$\sum_{k=0}^{p-1}(-1)^kA_k\eq\sum_{k=0}^{p-1}\f{\bi{2k}k^3}{16^k}\pmod{p^3}.\tag4.1$$
If $p\eq1,3\pmod8$, then
$$\sum_{k=0}^{p-1}A_k\eq\sum_{k=0}^{p-1}\f{\bi{4k}{k,k,k,k}}{256^k}\pmod{p^3}.\tag4.2$$

{\rm (ii)} If $x$ belongs to the set
$$\align &\{1,-4,9,-48,81,-324,2401,9801,-25920,-777924,96059601\}
\\&\quad\qquad\bigcup\l\{\f{81}{256},-\f9{16},\f{81}{32},-\f{3969}{256}\r\}
\endalign$$
and $x\not\eq0\pmod p$, then we must have
$$\sum_{k=0}^{p-1}A_k(x)\eq\l(\f xp\r)\sum_{k=0}^{p-1}\f{\bi{4k}{k,k,k,k}}{(256x)^k}\pmod{p^2}.$$
\endproclaim
\Remark\ 4.1. For
those
$$x=-4,9,-48,81,-324,2401,9801,-25920,-777924, 96059601, \f{81}{256},$$
the author (cf. [Su2])
had conjectures on $\sum_{k=0}^{p-1}\bi{4k}{k,k,k,k}/(256x)^k$ mod $p^2$. Motivated by this, Z. H. Sun [S2]
guessed
$\sum_{k=0}^{p-1}\bi{4k}{k,k,k,k}/(256x)^k$ mod $p^2$ for $x=-9/16,\,81/32,\,-3969/256$ in a similar way.

\medskip

We have checked Conjecture 4.1 as well as all the other conjectures in this paper via {\tt Mathematica 7}.
Below we provide numerical evidences for (4.1) and (4.2).
\smallskip

{\it Example} 4.1. The values of $A_0,A_1,\ldots,A_{10}$ are given by
$$\gather 1,\ 5,\ 73,\ 1445,\ 33001,\ 819005,\ 21460825,
\\584307365,\ 16367912425,\ 468690849005,\ 13657436403073\endgather$$
respectively. Via computation we find that
$$\sum_{k=0}^6(-1)^kA_k=20673445,\ \ \ \sum_{k=0}^6\f{\bi{2k}k^3}{16^k}=\f{18825543}{262144},$$
and also
$$20673445\times 262144-18825543=7^3\times 15800002159.$$
This verifies (4.1) for $p=7$. By computation we have
$$\sum_{k=0}^{10}A_k=14143101786223,\ \ \sum_{k=0}^{10}\f{\bi{4k}{k,k,k,k}}{256^k}
=\f{22821835381970859184405}{18889465931478580854784},$$
and also
$$\align&14143101786223\times 18889465931478580854784-22821835381970859184405
\\&\qquad=11^3\times 200717985992690007194123778899817.\endalign$$
This verifies (4.2) for $p=11$.

\medskip

Inspired by parts (ii) and (iii) of Theorem 1.1, we raise the following conjecture.

\proclaim{Conjecture 4.2} For any $\ve\in\{\pm1\}$, $m,n\in\Z^+$ and $x\in\Z$, we have
$$\sum_{k=0}^{n-1}(2k+1)\ve^k A_k(x)^m\eq0\ (\mo\ n).\tag4.3$$
If $p>5$ is a prime, then
$$\sum_{k=0}^{p-1}(2k+1)A_k\eq p-\f 72p^2H_{p-1}\pmod{p^6}.\tag4.4$$
\endproclaim

\Remark\ 4.2. After reading an initial version of this paper, Guo and Zeng [GZ] proved
the author's following conjectural results:

(a) For any $n\in\Z^+$ and $x\in\Z$ we have
$$\sum_{k=0}^{n-1}(2k+1)(-1)^k A_k(x)\eq0\ (\mo\ n).$$
If $p$ is an odd prime and $x$ is an integer, then
$$\sum_{k=0}^{p-1}(2k+1)(-1)^kA_k(x)\eq p\l(\f {1-4x}p\r)\ (\mo\ p^2).$$

(b) For any prime $p>3$ we have
$$\sum_{k=0}^{p-1}(2k+1)(-1)^kA_k\eq p\l(\f p3\r)\ (\mo\ p^3)$$
and
$$\sum_{k=0}^{p-1}(2k+1)(-1)^kA_k(-2)\eq p-\f 43p^2q_p(2)\ (\mo\ p^3).$$
\medskip

{\it Example} 4.2. It is easy to check that
$$\sum_{k=0}^9(2k+1)A_k^2=4178310699572329761604780\eq0\ (\mo\ 10)$$
and
$$\sum_{k=0}^9(2k+1)(-1)^kA_k^2=-4169201796654383947725970\eq0\ (\mo\ 10).$$
So (4.3) holds when $m=2,\ n=10,\ x=1$ and $\ve\in\{\pm1\}$.
Via computation we find that
$$\sum_{k=0}^{10}(2k+1)A_k=295998598024613,\ \ \ \ H_{10}=\f{7381}{2520},$$
and
$$295998598024613-11+\f 72\times 11^2\times\f{7381}{2520}=11^6\times\f{120300114181}{720}.$$
This verifies (4.4) for $p=11$.
\medskip

Recall that for a prime $p$ and a rational number $x$, the {\it $p$-adic
 valuation}
 of $x$ is given by
 $$\nu_p(x)=\sup\{a\in\Z:\ \t{the denominator of } p^{-a}x\ \t{is not divisible by}\ p\}.$$
Just like the Ap\'ery polynomial $A_n(x)=\sum_{k=0}^n\bi nk^2\bi{n+k}k^2x^k$ we define
$$D_n(x)=\sum_{k=0}^n\bi nk\bi{n+k}kx^k.$$
Actually $D_n((x-1)/2)$ coincides with the Legendre polynomial $P_n(x)$ of degree $n$.

Our following conjecture involves $p$-adic valuations.

\proclaim{Conjecture 4.3} {\rm (i)} For any $n\in\Z^+$ the numbers
$$s(n)=\f1{n^2}\sum_{k=0}^{n-1}(2k+1)(-1)^kA_k\l(\f14\r)$$
and
$$t(n)=\f1{n^2}\sum_{k=0}^{n-1}(2k+1)(-1)^kD_k\l(-\f14\r)^3$$
are rational numbers with denominators $2^{2\nu_2(n!)}$ and $2^{3(n-1+\nu_2(n!))-\nu_2(n)}$ respectively.
Moreover, the numerators of $s(1),s(3),s(5),\ldots$ are congruent to $1$ modulo $12$ and the numerators
of $s(2),s(4),s(6),\ldots$ are congruent to $7$ modulo $12$.
If $p$ is an odd prime and $a\in\Z^+$, then
$$s(p^a)\eq t(p^a)\eq1\ (\mo\ p).$$
For $p=3$ and $a\in\Z^+$ we have
$$s(3^a)\eq4\ (\mo\ 3^2)\ \ \ \t{and}\ \ \ t(3^a)\eq-8\ (\mo\ 3^5).$$

{\rm (ii)} Let $p$ be a prime. For any positive integer $n$ and $p$-adic integer $x$, we have
$$\nu_p\(\f1n\sum_{k=0}^{n-1}(2k+1)(-1)^kA_k\l(x\r)\)\gs\min\{\nu_p(n),\,\nu_p(4x-1)\}\tag4.5$$
and
$$\nu_p\(\f1n\sum_{k=0}^{n-1}(2k+1)(-1)^kD_k\l(x\r)^3\)\gs\min\{\nu_p(n),\,\nu_p(4x+1)\}.\tag4.6$$
\endproclaim

{\it Example}\ 4.3. We check Conjecture 4.3 with $n=p=5$. For $n=5$ we have
$$2^{2\nu_2(n!)}=2^{2\nu_2(120)}=2^6=64$$
and $$2^{3(n-1+\nu_2(n!))-\nu_2(n)}=2^{3(5-1+3)-0}=2^{21}=2097152.$$
Via computation we find that
$$s(5)=\f{19849}{64}\ \ \ \t{and}\ \ \ t(5)=\f{82547}{2097152}.$$
Note that $19849\eq1\ (\mo\ 12)$ and $s(5)\eq t(5)\eq1\ (\mo\ 5)$.
Also,
$$\nu_5\(\f15\sum_{k=0}^4(2k+1)(-1)^kA_k(-1)\)=\nu_5(-13095)=1=\nu_5(\pm5)$$
and
$$\nu_5\(\f15\sum_{k=0}^4(2k+1)(-1)^kD_k(1)^3\)=\nu_5(59189205)=1=\nu_5(5).$$
\medskip

Motivated by Theorem 1.3, we pose the following conjecture.

\proclaim{Conjecture 4.4} Let $p$ be an odd prime and let $n\gs 2$ be an integer.
Suppose that $x$ is a $p$-adic integer with $x\eq-2k\pmod{p}$ for some
$k\in\{1,\ldots,\lfloor(p+1)/(2n+1)\rfloor\}$. Then we have
$$\sum_{r=0}^{p-1}(-1)^r\bi xr^{2n+1}\eq0\pmod{p^2}.\tag4.7$$
\endproclaim

{\it Example} 4.4. Clearly $1/3\eq-2\ (\mo\ 7)$ and $\lfloor (7+1)/5\rfloor=1$. Via computation we find that
$$\sum_{r=0}^6(-1)^r\bi{1/3}r^5=\f{12107415300799972328}{12157665459056928801}=7^2\times\f{247090108179591272}{12157665459056928801}.$$
So (4.7) holds when $p=7$, $n=2$ and $x=1/3$.

\proclaim{Conjecture 4.5} Let $p\eq3\ (\mo\ 4)$ be a prime and let $m\in\Z$ with $p\nmid m(4m+1)$. Then
$$\sum_{k=0}^{p-1}\f{W_k(-m^2)}{(4m+1)^k}\eq\sum_{k=0}^{p-1}\f{A_k(-m^2/(4m+1))}{(4m+1)^k}
\eq\sum_{k=0}^{p-1}\f{\bi{2k}k^2}{(-16)^k}D_{2k}\l(\f1{4m}\r)\ (\mo\ p^2).$$
\endproclaim
\Ack. The author is grateful to the referee for helpful comments.


\widestnumber\key{BEW}

 \Refs

\ref\key A\by S. Ahlgren\paper Gaussian hypergeometric series and combinatorial congruences
\jour in: Symbolic computation, number theory, special functions, physics and combinatorics (Gainesville, FI, 1999),
pp. 1-12, Dev. Math., Vol. 4, Kluwer, Dordrecht, 2001\endref

\ref\key AO\by S. Ahlgren and K. Ono\paper A Gaussian hypergeometric series evaluation and Ap\'ery number congruences
\jour J. Reine Angew. Math.\vol 518\yr 2000\pages 187--212\endref

\ref\key Ap\by R. Ap\'ery\paper Irrationalit\'e de $\zeta(2)$ et $\zeta(3)$
\jour Ast\'erisque\vol 61\yr 1979\pages 11--13\endref

\ref\key B\by F. Beukers\paper Another congruence for the Ap\'ery numbers
\jour J. Number Theory\vol 25\yr 1987\pages 201--210\endref


\ref\key BEW\by B. C. Berndt, R. J. Evans and K. S. Williams
\book Gauss and Jacobi Sums\publ John Wiley \& Sons, 1998\endref

\ref\key CHV\by J.S. Caughman, C.R. Haithcock and J.J.P. Veerman
\paper A note on lattice chains and Delannoy numbers\jour Discrete Math.\vol 308\yr 2008\pages 2623--2628\endref


\ref\key CDE\by S. Chowla, B. Dwork and R. J. Evans\paper On the mod $p^2$ determination of $\bi{(p-1)/2}{(p-1)/4}$
\jour J. Number Theory\vol24\yr 1986\pages 188--196\endref


\ref\key G\by H. W. Gould\book Combinatorial Identities
\publ Morgantown Printing and Binding Co., 1972\endref

\ref\key GKP\by R. L. Graham, D. E. Knuth and O. Patashnik
 \book Concrete Mathematics\publ 2nd ed., Addison-Wesley, New York\yr 1994\endref

\ref\key GZ\by V. J. W. Guo and J. Zeng\paper Proof of some
conjectures of Z.-W. Sun on congruences for Ap\'ery polynomials\jour
J. Number Theory \vol 132\yr 2012\pages 1731--1740\endref

\ref\key I\by T. Ishikawa\paper Super congruence for the Ap\'ery numbers
\jour Nagoya Math. J.\vol 118\yr 1990\pages 195--202\endref

 \ref\key L\by E. Lehmer\paper On congruences involving Bernoulli numbers and the quotients
of Fermat and Wilson\jour Ann. of Math.\vol 39\yr 1938\pages 350--360\endref

\ref\key M03\by E. Mortenson\paper A supercongruence conjecture of Rodriguez-Villegas
for a certain truncated hypergeometric function
\jour J. Number Theory\vol 99\yr 2003\pages 139--147\endref


\ref\key M05\by E. Mortenson\paper Supercongruences for truncated  ${}_{n+1}\! F_n$
hypergeometric series with applications to certain weight three newforms
\jour Proc. Amer. Math. Soc.\vol 133\yr 2005\pages 321--330\endref


\ref\key O\by K. Ono\book Web of Modularity: Arithmetic of the Coefficients of Modular Forms and $q$-series
\publ Amer. Math. Soc., Providence, R.I., 2003\endref

\ref\key P\by H. Pan\paper On a generalization of Carlitz's congruence\jour Int. J. Mod. Math.
\vol 4\yr 2009\pages 87--93\endref

\ref\key PWZ\by M. Petkov\v sek, H. S. Wilf and D. Zeilberger\book $A=B$ \publ A K Peters, Wellesley, 1996\endref

\ref\key Po\by A. van der Poorten\paper A proof that Euler missed$\ldots$Ap\'ery's proof
of the irrationality of $\zeta(3)$
\jour Math. Intelligencer\vol 1\yr1978/79\pages 195--203\endref

\ref\key RV\by F. Rodriguez-Villegas\paper Hypergeometric families of Calabi-Yau manifolds
\jour {\rm in}: Calabi-Yau Varieties and Mirror Symmetry (Toronto, ON, 2001), pp. 223-231,
Fields Inst. Commun., {\bf 38}, Amer. Math. Soc., Providence, RI, 2003\endref

\ref\key S\by N. J. A. Sloane\paper {\rm Sequence A001850 in OEIS
(On-Line Encyclopedia of Integer Sequences)}
\jour {\tt http://oeis.org/A001850}\endref

\ref\key S1\by Z. H. Sun\paper  Congruences concerning Bernoulli numbers and Bernoulli polynomials
\jour Discrete Appl. Math.\vol 105\yr 2000\pages 193--223\endref

\ref\key S2\by Z. H. Sun\paper Congruences concerning Legendre
polynomials \jour Proc. Amer. Math. Soc. \vol 139\yr 2011\pages 1915--1929\endref

\ref\key S3\by Z. H. Sun\paper Congruences concerning Legendre
polynomials II \jour arXiv:1012.3898 \endref

\ref\key Su1\by Z. W. Sun\paper On congruences related to central binomial coefficients
\jour J. Number Theory \vol 131\yr 1011\pages 2219--2238\endref

\ref\key Su2\by Z. W. Sun\paper Super congruences and Euler numbers
\jour Sci. China Math.\vol 54\yr 2011\pages 2509--2535\endref

\ref\key Su3\by Z. W. Sun\paper On Delannoy numbers and Schroder numbers
\jour J. Number Theory\vol 131\yr 2011\pages 2387--2397\endref

\ref\key Su4\by Z. W. Sun\paper On sums involving products of three binomial coefficients
\jour preprint, arXiv:1012.3141\endref

\ref\key ST\by Z. W. Sun and R. Tauraso\paper New congruences for central binomial coefficients
\jour Adv. in Appl. Math.\vol 45\yr 2010\pages 125--148\endref


\endRefs
\enddocument